\documentclass[12pt]{article}

\usepackage{amssymb}

\usepackage{latexsym}

\begin{document}

%\parindent=0pt

%\parskip=6pt

%\newfont{\blb}{msbm10 scaled\magstep1}

\newtheorem{theorem}{Theorem}

\newtheorem{lemma}[theorem]{Lemma}

\newtheorem{corollary}[theorem]{Corollary}

\newtheorem{proposition}[theorem]{Proposition}

\newtheorem{conjecture}[theorem]{Conjecture}

\def\deg{{\rm deg}}

\def\e{\epsilon}

\def\ve{\varepsilon}

\def\cA{{\cal A}}

\def\cG{{\cal G}}

\def\cS{{\cal S}}

\def\cL{{\cal L}}

\def\cK{{\cal K}}

\def\cF{{\cal F}}

\def\cP{{\cal P}}

\def\cI{{\cal I}}

\def\ex{{\rm ex}}

\def\cH{{\cal H}}

\def\cT{{\cal T}}

\title{\bf Finding bipartite subgraphs efficiently}

\author{Dhruv Mubayi\thanks{Department of Mathematics, Statistics, and Computer
Science, University of Illinois at Chicago, IL 60607;  email:
mubayi@math.uic.edu; research  supported in part by  NSF grant DMS
0653946}\quad  and \quad Gy\"orgy Tur\'an\thanks{ Department of
Mathematics, Statistics, and Computer Science, University of
Illinois at Chicago, IL 60607 and Research Group on Artificial
Intelligence, Hungarian Academy of Sciences, University of Szeged;
email:
 gyt@uic.edu}}
\date{\today}

%{Dhruv Mubayi and Gy\"orgy Tur\'an \thanks{

%Department of Mathematics, Statistics, \& Computer Science,

%University of Illinois at Chicago, 851 S. Morgan Street, Chicago, IL

%60607-7045; email: mubayi@math.uic.edu. \hfil\break\null\hskip .23in

%Keywords: {\it  }}  }

\date{}

%\today

\maketitle

\smallskip

\begin{abstract}

Polynomial algorithms are given for the following two problems:

\begin{itemize}
\item given a graph with $n$ vertices and $m$ edges, where $m
\ge 3 n^{3/2}$, find a complete balanced bipartite subgraph with
parts about $\frac{\ln n}{\ln (n^2/m)}$,

\item given a graph with $n$ vertices, find a decomposition of its edges
into complete balanced bipartite graphs having altogether $O(n^2 /
\ln n)$ vertices.
\end{itemize}

Previous proofs of the existence of such objects, due to K\H{o}v\'ari-S\'os-Tur\'an \cite{kst}, Chung-Erd\H{o}s-Spencer \cite{chu}, Bublitz \cite{bub} and Tuza \cite{tuz} were
non-constructive.
\end{abstract}

\maketitle

\section{Introduction}

Determining the minimal number of edges in a bipartite graph which
guarantees the existence of a complete balanced bipartite subgraph
$K_{q, q}$ is known as the Zarankiewicz problem (see, e.g.,
Bollob\'as \cite{boll}). It was shown by K\H{o}v\'ari, S\'os and
Tur\'an \cite{kst} that every bipartite graph with $n$ vertices in
both sides and $c_q n^{2 - 1/q}$ edges contains a $K_{q,q}$. The
same bound (with different constant $c_q$)  holds for general $n$-vertex graphs. The argument from
\cite{kst} also shows that $n$-vertex graphs of constant density,
i.e., graphs with $\e n^2$ edges, contain a complete bipartite
graph with parts of size at least $c_\e \ln n$. The proofs of all
these results are based on counting, and thus are
\emph{non-constructive}.

We consider the question whether such subgraphs can be found by
\emph{efficient}, i.e., polynomial time, algorithms. This question
has been considered recently by Kirchner \cite{kir}, who gave an
efficient algorithm to find a complete balanced bipartite subgraph
with parts of size $\Omega(\sqrt{\ln n})$ in graphs of constant
density. We improve this result by giving an efficient algorithm
which finds a complete balanced bipartite subgraph with parts of
size $\Omega(\ln n)$, i.e. of the optimal order of magnitude, in
graphs of constant density. Our algorithm gives subgraphs of
similar size as the counting argument in other ranges as well
\footnote{Note that the problem becomes meaningless in the sense
studied here for fewer than $n^{3/2}$ edges, as such graphs do not
always contain even $K_{2,2}$ subgraphs.}.

Finding a largest balanced complete bipartite subgraph is an
important optimization problem, which is known to be NP-hard, and
even hard to approximate (see, e.g., Feige and Kogan \cite{fei}).
We would like to emphasize that we are \emph{not} trying to give
an approximation algorithm for this problem. Our objective is to
give an efficient algorithm which finds a balanced complete
bipartite subgraph of size close to the largest size that is
guaranteed to exist knowing \emph{only} the number of edges in the
graph. Thus, even in a dense graph, we are finding a subgraph of
logarithmic size only. Results of this type are given, for
example, in Alon \emph{et al.} \cite{alon}.

The counting argument of \cite{kst} has several applications to
other combinatorial problems. It seems to be an interesting
question whether the algorithmic version of the counting argument
leads to further algorithmic results in these applications. As a
case in point, we consider the question of decomposing, or
partitioning, the edge set of a graph into complete bipartite
graphs. The motivation to look for such algorithms comes from an
application in approximation algorithms \cite{bha}.

Every $n$-vertex graph can be decomposed into at most $n - 1$
stars, and Graham and Pollak \cite{gra} showed that $n - 1$
complete bipartite graphs are necessary for the $n$-vertex
complete graph. Instead of minimizing the number of complete
bipartite graphs in a decomposition, one can also try to minimize
the complexity of decompositions, measured by \emph{the sum of the
number of vertices of the complete bipartite graphs used in the
decomposition}. This measure of complexity was suggested by
Tarj\'an \cite{tar} in the context of circuit complexity. For
recent connections to circuit complexity see Jukna \cite{juk}.

It was shown by Chung, Erd\H{o}s and Spencer \cite{chu}, and by
Bublitz \cite{bub}, that there is always a decomposition of
complexity $O(n^2/\ln n)$, and this order of magnitude is best
possible. Similar results were obtained by Tuza \cite{tuz} for
decomposing bipartite graphs.  These results are obtained by
repeatedly applying the counting argument to show the existence of
a large complete bipartite graph and removing its edges. Thus the
decomposition results obtained in \cite{bub,chu,tuz} are also
non-constructive. As a direct application of our algorithm for
finding bipartite subgraphs, we obtain efficient algorithms to
find decompositions of complexity $O(n^2/\ln n)$.

\section{Complete balanced bipartite subgraphs}

Searching for a $K_{q,q}$ by checking all subgraphs of that size
would give an algorithm with superpolynomial running time if $q$
is, say, logarithmic in the number of vertices. A polynomial
algorithm could be given by restricting the search space to a
polynomial size set of candidate subgraphs. One possibility for
that would be to find a bipartite subgraph $(R, S)$ with the
following properties:

\begin{itemize}
\item it is dense enough for the known results to guarantee the
\emph{existence} of a $K_{q,q}$, and

\item the number of $q$-element subsets of $R$ is only
\emph{polynomial}.
\end{itemize}

If such an $(R, S)$ can be found efficiently then a required
$K_{q,q}$ is obtained by checking the common neighborhood of all
$q$-element subsets of $R$. It turns out that this approach indeed
works if one chooses $R$ to be the right number of vertices with
maximal degree and $S$ to be the remaining vertices. Thus, we
consider the following algorithm, where $q(n, m)$ and $r(n, m)$
are functions to be determined.

\bigskip

Algorithm {\bf FIND-BIPARTITE}

\bigskip

{\bf input:} $G = (V, E)$ with $|V|=n$ and $|E|=m$

\bigskip

%$R := \emptyset$

%{\bf for} $i := 1$ {\bf to} $r$ {\bf do}

%\hspace{1cm} let $v$ be a vertex such that $|N(v) - R|$ is maximal

%\hspace{1cm} $R := R \cup \{v\}$

$q := q(n, m), \,\, r := r(n, m)$,

\medskip

$R := \,\, \textrm{$r$ vertices having highest degree}$

\medskip

{\bf for} all subsets $C \subseteq R$ with $|C| = q$ {\bf do}

\medskip

\hspace{0.5cm} $D := \bigcap \{N(v) - R \, : \, v \in C\}$

\medskip

\hspace{0.5cm} {\bf if} $|D| \ge q$ {\bf then} $D' := \,\,
\textrm{the first $q$ elements of $D$}$, \,\, {\bf return} $(C,
D')$

\bigskip

We now show that with the appropriate choice of $q(n, m)$ and
$r(n, m)$ the algorithm works.

\begin{theorem} \label{thm:subg}
Let

\[q := \left\lfloor \frac{\ln (n/2)}{\ln (2 e
n^2/m)}\right\rfloor, \,\,\,\, r := \left\lfloor \frac{q
n^2}{m}\right\rfloor. \]

If $n$ is sufficiently large and $m \ge 3 n^{3/2}$ then Algorithm
{\bf FIND-BIPARTITE} returns a $K_{q,q}$ (with $q \ge 2$ as long
as $m> 8 n^{3/2}$). The running time of the algorithm is
polynomial in $n$.
\end{theorem}

{\bf Remark.}  Note that our algorithm finds a $K_{q,q}$ in an $n$-vertex graph
with $m=c_qn^{2-1/q}$ edges as long as $c_q$ is large. This is optimal for $q=2,3$ as there exist $n$-vertex graphs with $c'_qn^{2-1/q}$ edges and no $K_{q,q}$, and if certain conjectures in extremal graph theory are true, then it is also optimal for $q>3$.

\noindent
{\bf Proof.} After selecting $i < r$ vertices, the number of edges
incident to these vertices is less than $r n$. Hence in the
subgraph induced by the remaining vertices there is a vertex of
degree at least $2 (m - r n) / n$. Thus if $R$ is the set of $r$
highest degree vertices in $G$ then
\[ \sum_{v \in R} deg_G(v) \ge \frac{2 r (m - r n)}{n}. \]

%Let $R_i$ be the set $R$ after $i$ iterations, and $v_i$ be the
%vertex selected in the $i$'th iteration. Then the number of edges
%with no endpoint in $R_i$ is at least $m - i n \ge m - r n$, and
%so
%\[ N(v_{i}) - R_{i-1} \ge \frac{m - r n}{n}. \]
%In subsequent iterations the number of neighbors of $v_i$ outside
%$R$ may decrease, but at the end it still holds that
%\[ N(v_{i}) - R \ge \frac{m - 2 r n}{n} \ge \frac{m}{2 n}, \]
%where the last inequality follows from the assumption $m \ge 3
%n^{3/2}$.

Hence the bipartite graph $H$ with parts $R, V-R$ and edge set
comprising those edges of $G$ with one endpoint in $R$ and the
other in $V-R$ has at least $2 r m / n - 3r^2$ edges.

We will now argue that
$rm/n\ge 3r^2$.
Indeed, $rm/n\ge 3r^2$ is equivalent to $r\le m/3n$.
Now $r \le qn^2/m$ so it is enough to show that $qn^2/m \le m/3n$ or equivalently, that
$3qn^3\le m^2$.  Using the definition of $q$, we see that $3qn^3\le m^2$ follows from
$$m^2\ln (2en^2/m)\ge 3n^3\ln (n/2).$$  Suppose first that $3 n^{3/2} \le m \le 3
n^{3/2} \sqrt{\ln n}$.  Then
$$m^2\ln\left(\frac{2en^2}{m}\right)\ge 9n^3\ln\left(\frac{2en^2}{3n^{3/2}\sqrt{\ln n}}\right)> 9n^3\ln\left(\sqrt{\frac{n}{\ln n}}\right)>4n^3\ln n >3n^3\ln(n/2).$$
On the other hand, if
 $m \ge 3 n^{3/2} \sqrt{\ln n}$, then using $m<n^2/2$ we have
 $$m^2\ln(2en^2/m)\ge 9n^3\ln n\ln(2en^2/m)>
  9n^3\ln n \ln(4e)>3n^3\ln (n/2).$$
We conclude that $H$ has at least $2rm/n-3r^2\ge rm / n$ edges.

For the correctness of the algorithm it is sufficient to show that
$H$ contains a copy of $K_{q,q}$. This follows by the counting
argument referred to in the introduction, which is included here
for completeness. Let $s$ denote the number of stars with centers
in $V-R$ and $q$ leaves. Then

\[ s = \sum_{v \in V-R} {deg_H(v)\choose q} \,\, \ge \,\,
\frac{n}{2}{rm/n^2\choose q}, \] using the convexity of the
function which is ${x \choose q}$ if $x \ge q - 1$ and 0
otherwise, and using $r \le n/2$ which follows by the lower bound
on $m$. If the latter quantity is greater than $(q - 1) {r \choose
q}$ then there is a $q$-subset of $R$ which is the leaf set for at
least $q$ distinct stars, and this gives a copy of $K_{q,q}$.  Observe that the definition of $q$ implies that $n/2\ge (2en^2/m)^q$ and this is equivalent to
$$\frac{n}{2} \left(\frac{r m}{n^2 q}\right)^q \ge
\left(\frac{2 e r}{q}\right)^q.$$
Now the inequality above and standard estimates of the binomial coefficients give
\[ \frac{n}{2}{rm/n^2\choose q} > \frac{n}{2} \left(\frac{r m}{n^2 q}\right)^q \ge
\left(\frac{2 e r}{q}\right)^q \ge q \left(\frac{r e}{q}\right)^q >
(q - 1) {r \choose q},
\] Thus $H$ indeed contains a $K_{q,q}$.
%\[>\frac{n}{2}\left(\frac{m}{8n^2}\right)^q\left(\frac{re}{q}\right)^q>\frac{n}{2}\left(\frac{m}{8n^2}\right)^q{r
%\choose q}>q{r \choose q}\]

In order to show that the running time of the algorithm is
polynomial, note first that, assuming an adjacency matrix
representation, the set $R$ can be found in $O(n^2)$ steps. For a
given $q$-subset of $R$, the common neighbors can be found in $O(n
q)$ steps. All $q$-subsets can be listed in $O({r \choose q})$
steps (see, e.g. \cite{rein}). Thus the algorithm requires time
\[ O\left(n^2 + {r \choose q} n q\right). \]
The number of iterations is at most
\[ {r \choose q} \le \left(\frac{r e}{q}\right)^q \le e^q \left(\frac{n^2}{m}\right)^q =
e^q e^{q\ln(n^2/m)}\]
Now $m<n^2/2$ implies that
$$e^q\le e^{\ln n/\ln 4e}=n^{1/\ln 4e}<n^{0.4195}.$$
and $q<\ln n/\ln(n^2/m)$ implies that
$$e^{q \ln(n^2/m)}< e^{\ln n}=n.$$
Therefore the
running time of the algorithm is $O(n^{2.42})$. $\Box$

\section{Decomposition into balanced complete bipartite subgraphs}

Given a graph $G = (V, E)$, we consider complete bipartite
subgraphs $G_i = (A_i, B_i, E_i), i = 1, \ldots, t$ such that the
edges sets $E_i$ form a partition of $E$. The \emph{complexity} of
such a decomposition is measured by the total number of vertices,
i.e., by
\[ \sum_{i = 1}^t |A_i| + |B_i|. \]
We find a decomposition of complexity $O(n^2/\ln n)$. The
decomposition contains \emph{balanced} bipartite graphs, thus
$|A_i| = |B_i|$ holds as well. The algorithm uses Algorithm {\bf
FIND-BIPARTITE} in a straightforward manner. As stated, Algorithm
{\bf FIND-BIPARTITE} is guaranteed to work only if $n \ge n_0$ for
some $n_0$. As we are only interested in proving an asymptotic
result, let us assume that graphs on fewer vertices are handled by
some brute-force method.

\bigskip

Algorithm {\bf FIND-DECOMPOSITION}

\bigskip

\begin{quote}

Given an $n$-vertex input graph $G = (V, E)$, if $n < n_0$, use a
brute-force method to find an optimal decomposition of $G$. Else,
use Algorithm {\bf FIND-BIPARTITE} repeatedly to find a complete
balanced bipartite subgraph and delete it from the current graph,
as long as there are more than $n^2 / \ln n$ edges. After that,
form a separate bipartite graph from each remaining edge.
\end{quote}

\bigskip

\begin{theorem} \label{thm:deco}
For every $n$-vertex graph $G$, Algorithm {\bf FIND-DECOMPOSITION}
finds a decomposition of $G$ into balanced complete bipartite
graphs, having complexity
\[ O\left(\frac{n^2}{\ln n}\right). \]
 The running time of the algorithm is polynomial in $n$.
\end{theorem}

\noindent
{\bf Proof.} As the size of the subgraphs produced by Algorithm
{\bf FIND-BIPARTITE} is of the same order of magnitude as
guaranteed by the existence theorems, the theorem follows as in
\cite{bub,chu,tuz}. For completeness, we give the argument,
following \cite{tuz}.

Let the subgraphs produced by the calls of  Algorithm {\bf
FIND-BIPARTITE} be $G_i = (A_i, B_i)$ with $|A_i| = |B_i| = q_i$,
where $i = 1, \ldots, t$ for some $t$ . We need to show that
\begin{equation} \sum_i q_i = O\left(\frac{n^2}{\ln n}\right). \label{eq:charg}
\end{equation}
Let us divide the iterations of the algorithm into \emph{phases}.
The $\ell$'th phase consists of those iterations where the number
of edges in the input graph of Algorithm  {\bf FIND-BIPARTITE} is
more than ${n^2}/ (\ell + 1)$ and at most ${n^2}/ \ell$. Dividing
up the term $q_i$ in (\ref{eq:charg}) between the $q_i^2$ edges of
$G_i$, each edge gets a weight of $1 / q_i$. We have to upper
bound the sum of the weights assigned to the edges.

It follows from the definition of $q_i$ in Theorem \ref{thm:subg}
that graphs formed in the $\ell$'th phase have $q_i = \Theta(\ln
n / \ln \ell)$. Thus edges, which get their weight in the
$\ell$'th phase, get a weight of $\Theta(\ln \ell/ \ln n)$. The
number of edges getting their weight in the $\ell$'th phase is
$\Theta((\frac{1}{\ell} - \frac{1}{\ell + 1}) n^2) = \Theta(n^2 /
\ell^2)$. Hence the total weight assigned to the edges is at most
of the order of magnitude
\[ \sum_{\ell = 1}^\infty \frac{\ln \ell}{\ln n} \cdot \frac{n^2}{\ell^2}
 = \Theta\left(\frac{n^2}{\ln n}\right), \] as $\sum
\frac{\ln \ell}{\ell^2}$ is convergent. The polynomiality of the
running time follows directly from the polynomial running time of
Algorithm  {\bf FIND-BIPARTITE}. $\Box$

\section{Subgraphs and decompositions of bipartite graphs}

In this section we formulate the result analogous to Theorem
\ref{thm:deco} for bipartite graphs $G = (A, B,
E)$ having parts $A$ and $B$, with $|A| = a$, $|B| = b$ and $|E| =
m$. We assume \emph{w.l.o.g.} that $a \ge b$.

The algorithms and their analysis are straightforward
modifications of those for general graphs. Algorithm {\bf
FIND-BIPARTITE-IN-BIPARTITE}, a modified version of {\bf
FIND-BIPARTITE}, uses functions $q(a, b, m)$ and $r(a, b, m)$. It
constructs $R$ as the set of $r$ highest degree vertices in $B$,
and checks the common neighborhood of all $q$ element subsets of
$R$. Algorithm {\bf FIND-DECOMPOSITION-IN-BIPARTITE}, a modified
version of {\bf FIND-DECOMPOSITION}, uses this modified algorithm
while the number of edges is greater than $a b / \ln (a+b)$.

\begin{theorem}
Let $G$ be a bipartite graph with sides of size $a$ and $b$.
 Algorithm  {\bf FIND-DECOMPOSITION-IN-BIPARTITE} finds a
decomposition of $G$ into balanced complete bipartite graphs,
having complexity
\[ O\left(\frac{a b}{\ln (a+b)}\right). \]
 The running time of the algorithm is polynomial in $a + b$.
\end{theorem}

\bigskip

{\bf Acknowledgment} We thank Stefan Kirchner for sending us his
Ph.D. dissertation.


\begin{thebibliography}{99}
\bibitem{alon} N. Alon, R. A. Duke, H. Lefmann, V. R\"odl, R.
Yuster: The algorithmic aspects of the regularity lemma, \emph{J.
of Algorithms} {\bf 16} (1994), 80-109.

\bibitem{bha} A. Bhattacharya, B. DasGupta, Gy. Tur\'an: On
approximate Horn minimization. In preparation.

\bibitem{boll} B. Bollob\'as: \emph{Extremal Graph Theory}.
Academic Press, 1978.

\bibitem{bub} S. Bublitz: Decomposition of graphs and monotone
formula size of homogeneous fuctions, \emph{Acta Informatica} {\bf
23} (1986), 689-696.

\bibitem{chu} F. R. K. Chung, P. Erd\H{o}s, J. Spencer: On the
decomposition of graphs into complete bipartite graphs, in:
\emph{Studies in Pure Mathematics, To the Memory of Paul Tur\'an},
95-101. Akad\'emiai Kiad\'o, 1983.

\bibitem{fei} U. Feige, S. Kogan: Hardness of approximation of the
balanced complete bipartite subgraph problem, \emph{Tech. Rep.
MCS04-04, Dept. of Comp. Sci. and Appl. Math., The Weizmann Inst.
of Science}, 2004.

\bibitem{gra} R. L. Graham, H. O. Pollak: On the addressing
problem for loop switching, \emph{Bell Syst. Techn. J.} {\bf 50}
(1971), 2495-2519.

\bibitem{juk} S. Jukna: Disproving the single level conjecture,
\emph{SIAM J. Comp.} {\bf 36} (2006), 83-98.

\bibitem{kir} S. Kirchner: Lower bounds for Steiner tree
algorithms and the construction of bicliques in dense graphs.
Ph.D. Dissertation, Humboldt-Universit\"at zu Berlin, 2008. (In
German.)

\bibitem{kst} T. K\H{o}v\'ari, V. T. S\'os, P. Tur\'an: On a
problem of K. Zarankiewicz, \emph{Colloq. Math.} {\bf 3} (1954),
50-57.

\bibitem{rein} E. M. Reingold, J. Nievergelt, N. Deo:
\emph{Combinatorial Algorithms}. Prentice Hall, 1977.


\bibitem{tar} T. Tarj\'an: Complexity of lattice-configurations,
\emph{Studia Sci. Math. Hung.} {\bf 10} (1975), 203-211.

\bibitem{tuz} Zs. Tuza: Covering of graphs by complete bipartite
subgraphs; complexity of 0-1 matrices, \emph{Combinatorica} {\bf
4} (1984), 111-116.


\end{thebibliography}
\end{document}